\documentclass[11pt]{amsart}

\usepackage{hyperref,wasysym}
\hypersetup{colorlinks=true, citecolor=BurntOrange, linkcolor=ForestGreen}

\usepackage{amssymb}
\usepackage{graphicx}
\usepackage{amsmath}
\usepackage[usenames,dvipsnames]{xcolor}
\usepackage[text={6.5in,9in},centering]{geometry}
\newtheorem{Lemma}{Lemma}
\newtheorem{Theorem}[Lemma]{Theorem}

\newenvironment{Proof}[1][.]%
  {\begin{trivlist}\item[]\textbf{Proof#1 }}%
  {\hspace*{\fill}$\rule{.3\baselineskip}{.35\baselineskip}$\end{trivlist}}

\newcommand{\eps}{\epsilon} 

\begin{document}

\title{ Fisher-KPP dynamics in diffusive Rosenzweig-MacArthur and Holling-Tanner  models.}

%  equation Traveling fronts in a modified diffusive  Rosenzweig-MacArthur system}

%\author[Ghazaryan]{
 
\date{\today}

\begin{abstract}{ }
We prove the existence of traveling fronts in  diffusive Rosenzweig-MacArthur   and Holling-Tanner population models and investigate their relation with fronts in a scalar Fisher-KPP equation. 
More precisely, we prove the existence of fronts  in a   Rosenzweig-MacArthur predator-prey model in  two situations: when the prey diffuses at the rate much smaller than that of the predator and when  both the predator and the prey diffuse very slowly. Both situations are captured as singular perturbations of  the associated limiting systems. In the first situation we demonstrate  clear relations of the fronts with the fronts in a scalar Fisher-KPP equation. 
Indeed, we show that  the underlying dynamical system in  a singular limit is reduced to a scalar Fisher-KPP equation and the fronts supported by the full  system are small perturbations of the Fisher-KPP fronts.
We obtain a  similar result for a diffusive Holling-Tanner population model. In the second situation for the  Rosenzweig-MacArthur model we prove the existence of the fronts but  without observing a direct relation with Fisher-KPP equation. The analysis suggests that, in a variety of  reaction-diffusion  systems that rise in population modeling,    parameter regimes may be found when the   dynamics of the system  is inherited from  the  scalar Fisher-KPP equation.
\end{abstract}
\maketitle

\begin{center} Hong Cai  \textsuperscript{a},  Anna Ghazaryan \textsuperscript{b},  Vahagn Manukian \textsuperscript{b,c}\\
 \end{center}
 
  \textsuperscript{a} \address{ %Hong_Cai@brown.edu 
  Department of Physics, Brown University, 182 Hope Street, Providence, RI 02912, USA,
  %Department of Mathematics, Miami University, 301 S. Patterson Ave,   Oxford, OH 45056, USA,
  }  
 \email{Hong\_Cai@brown.edu }\par
 \textsuperscript{b} \address{ Department of Mathematics, Miami University, 301 S. Patterson Ave,   Oxford, OH 45056, USA,}  
 \email{ghazarar@miamioh.edu} \par
 \textsuperscript{c} \address{ Department of Mathematical and Physical Sciences, Miami University, 1601 University Blvd, Hamilton, OH 45011, USA,}  \email{manukive@miamioh.edu}


%35B40 %Asymptotic behavior of solutions,
 
\section{Plan of the paper}

Reaction-diffusion systems are often used in population dynamics modeling when it is desirable to take into account random motion of individuals in the population. In systems where  there is more than one spatially homogeneous equilibrium state, it is of interest to know whether  transition fronts between these states exist.  
We present analysis of traveling fronts in two diffusive population-dynamics models, one of which is a modified Rosenzweig-MacArthur model,  and the other is Holling-Tanner predator-prey model. The analysis is performed in detail on  the  Rosenzweig-MacArthur  system (Section~\ref{s:intro}), while  for the Holling-Tanner system (Section~\ref{HT}) most of the details are skipped and similarities in the proofs are pointed out.  The plan of the paper is as follows. We introduce  the Rosenzweig-MacArthur model in Section~\ref{background}  and explain  the results of the paper and their mathematical implications. The scaling of the model that we use and the parameter regimes  that the analysis covers  are described in Section~\ref{RM}. The regimes are grouped in  two cases that are then analyzed using geometric singular perturbation theory in Sections~\ref{sec:1d} and \ref{sec:2d}. In Section~\ref{sec:1d} the relation of the fronts  with a scalar Fisher-KPP equation is revealed. In the analysis of the Rosenzweig-MacArthur model, assumptions are made about some of the parameters, in order to simplify some of the calculations.   In Section~\ref{gamma}, we describe the implications of the assumptions and how the obtained results extend  to the complementary cases.
We  then, in Section~\ref{HT}, describe the  Holling-Tanner population model and  a regime  where  the fronts are constructed in a  way similar to the one implemented in Section \ref{sec:1d}.

\section{Rosenzweig-MacArthur model  \label{s:intro}}
\subsection{ The  background and  physical interpretation\label{background}}

In 1963, Rosenzweig and MacArthur  suggested  \cite{RM} that  in predator-prey interactions the predator's ability to increase has a ``ceiling", i.e. there is a bound on the rate at which predators consume  the prey,  regardless of the density of the prey population.    Population models that take this fact into account are called  Rosenzweig-MacArthur  models. 
A classical Rosenzweig-MacArthur model  may have a form of a  system of ordinary differential equations such as
\begin{eqnarray}\label{e:ode0}
U_{\tau}&=&\mathcal AU\left(1-\frac{U}{\mathcal K}\right)-\frac{\mathcal B UW}{\mathcal E+U},\notag\\
W_{\tau}&=&- CW+ \frac{\mathcal D UW}{\mathcal E+  U},
\end{eqnarray} 
where $\tau$ is the time,  and $U$ and $W$ are the population densities of the prey and predator, respectively. 
 Parameter $\mathcal A>0$ is the growth factor for the prey species, $\mathcal C>0$ is the death rate for the predator without prey, $\mathcal K>0$ is the carrying capacity of the prey species. Positive parameters $\mathcal B$ and $\mathcal D$ are the interaction rates for the two species. The term $\frac{UW}{\mathcal E+U}$, where $\mathcal E>0$ is a constant, represents the MacArthur-Rosenzweig effect.
 
In real life, mechanisms of the phenomenon that Rosenzweig-MacArthur captures mathematically  may be various and complex. A somewhat simplified physical interpretation of the Rosenzweig-MacArthur   effect  for the prey is  that sometimes higher predation rates  cause or are associated with a decrease in the prey mortality.  One  example is  related to   the ability of the prey to   take an environmental refuge. The  effect  of the environmental protection  of the prey on the predator is  negative.

Variations of \eqref{e:ode0} exist. Although the Rosenzweig-MacArthur  system  is at this point a classical  predator-prey model,  this system and its variations still attract attention  of mathematicians as it supports a plethora of interesting mathematical phenomena. For example, periodic wave trains have been numerically observed in \cite{SS1}.  Another example  is a paper  \cite{AVOOMR} that contains  a study of  a Rosenzweig-MacArthur  model with a ``refuge function" proposed by Almanza-Vasquez  in \cite{AVGOGY}. More precisely,  it is suggested that  if a quantity $U_r$ of prey population takes a refuge, then only  $U-U_r$ prey  is exposed to the predator,  and so  the  term $\frac{ UW}{\mathcal E+U} $  in \eqref{e:ode0}  is replaced by  a term of the form $\frac{ (U-U_r) W}{\mathcal E+(U-U_r)}$. In this context of the refuge, $\mathcal E$ is the number of prey  that represents a half of the maximum capacity of the refuge.  The  paper contains a study  of the  stability of  the physically meaningful equilibrium points, the influence of the size of the refuge on  the coexistence equilibrium  is observed and, moreover, the  existence of limit cycles resulting in oscillations in populations of predator and prey is shown. 

The concept of the refuge is not the only phenomenon captured mathematically by Rosenzweig-MacArthur systems.   A different  example  is formulated in \cite{AADO}:
 \begin{eqnarray}\label{e:010}
U_{\tau}&=&\mathcal AU\left(1-\frac{U}{\mathcal K}\right)-\frac{\mathcal B UW}{1+\mathcal E_1U},\notag\\
W_{\tau}&=&-\mathcal CW+ \frac{\mathcal D UW}{1+ \mathcal E_2U}.
\end{eqnarray} 
Here, the  two  positive  constants $\mathcal E_1$ and $\mathcal E_2$ are related to how much  protection the environment provides to the  prey  and predator, respectively. Paper  \cite{AADO}   is devoted to the study of the
 the boundedness  of the solutions and global stability  of  the co-existence equilibrium for a predator-prey model   \cite{AADO}.
 According to \cite{AADO}, the  model may represent,  for example,  an insect pest-spider food chain. Other   examples  of   real-life  population systems may be found in \cite{6,SGO}.

The  systems of  local  equations  such as Rosenzweig-MacArthur  model provide a valuable insight   into the population dynamics. 
It is also widely accepted that  reaction-diffusion systems or partly parabolic systems (systems where some but  not all  quantities diffuse) built as extensions of local  predator-prey interaction  models and other biological models  are of great interest as well (see for example \cite{VP} and the references therein).  In reaction-diffusion population models one would be interested in establishing the existence of  traveling wave solutions such as periodic wavetrains, pulses, or fronts.  We refer readers to a review paper \cite{SS} and references therein on periodic traveling waves in cyclic populations. For fronts,  we mention the  paper \cite{Dunbar}  where Dunbar considered  the system
\begin{eqnarray} \label{dunbar}
U_{\tau}&=&\epsilon_u U_{XX}+\mathcal AU\left(1-\frac{U}{\mathcal K}\right)-\frac{\mathcal B UW}{1+ \mathcal EU},\notag\\
W_{\tau}&=&\epsilon_wW_{XX}-\mathcal CW+ \frac{\mathcal D UW}{ 1+\mathcal E U}.
\end{eqnarray} 
Here $X$ is a one-dimensional spatial variable  and  $\mathcal E$  captures the  following satiation phenomenon: ``the consumption of prey by a unit number of predators cannot continue to grow linearly with the number of prey available but must ``saturate" at the value $1/\mathcal E$". For the case $\epsilon_u=0$, Dunbar shows  the existence of the periodic traveling wavetrains and  fronts connecting an equilibrium to a periodic orbit  using shooting techniques, invariant manifold theory, and the qualitative theory of ordinary differential equations.

The system \eqref{dunbar} has been well studied:  In \cite{OL} in  the case of non-zero $\epsilon_u$ and $\epsilon_w$  the existence of traveling waves was shown numerically. In \cite{HLR}   the existence of traveling wave solutions and small amplitude traveling wave train solutions of system  was proved using  a  shooting argument together with   a Liapunov function and LaSalle's Invariance Principle.

Dunbar in \cite{Dunbar} mentions that %it is not important that equations are  rooted in mathematical ecology since 
analogous  equations appear in  a variety  disciplines   such as chemical kinetics, cell biology and immunology.  In particular,  the  term that represents the satiation effect  appears in cell biology (Michaelis-Menten dynamics).  According to Dunbar, the results must be perceived beyond the implementations to a particular model and  are ``meant to indicate some of the interesting ranges of behavior possible for systems of reaction-diffusion equations".

The  interest in learning whether  traveling wave solutions  in biological systems exist or not dates back to the important works of   Fisher \cite{Fisher} and Kolmogorov et al. \cite{KPP} on Fisher-KPP equation. Traveling waves in biological systems are cornerstones of the classic books  by  Fife \cite{Fife}, Murray \cite{M},  and Volpert et al. \cite{VVV}  that contain numerous references. 
In  particular, traveling fronts in  biological reaction-diffusion systems  are of interest. Introducing a diffusion term  in a reaction equation often allows to capture a moving zone of transition between an  absence of  a population to a nonzero equilibrium state.

In the current paper, we consider a diffusive Rosenzweig-MacArthur system that we describe in the next section. For this system, we give a proof of the existence of traveling fronts in certain parameter regimes by  geometric construction. Through this geometric construction we demonstrate  in Section 3 that  in some  parameter regimes the dynamics of the system is driven by a scalar Fisher-KPP  equation. On the other hand,  in a parameter regime considered in Section 4, we do not observe the relation of the  fronts  with the Fisher-KPP equation. % This makes the implications of this analysis even more interesting. For example, 
Comparing the stability properties  of these two different classes of fronts is a work in progress.
% we suggest as  an open problem. 
  We also mention that transition fronts in the Rosenzweig-MacArthur system have been observed numerically and their stability properties have been investigated  in \cite{DS}, but in parameter regimes which are not covered in the current paper.

\subsection{The model\label{RM}}

In the current paper, we consider a version of a diffusive Rosenzweig-MacArthur system  that is based on the  reaction system considered in \cite{AADO}, 
\begin{eqnarray}\label{e:10}
U_{\tau}&=&\epsilon_u U_{XX}+\mathcal AU\left(1-\frac{U}{\mathcal K}\right)-\frac{\mathcal B UW}{1+\mathcal E_1U},\notag\\
W_{\tau}&=&\epsilon_wW_{XX}-\mathcal CW+ \frac{\mathcal D UW}{1+ \mathcal E_2U}.
\end{eqnarray} 
We reiterate that, here,  $\tau\ge1 0$ is the time, $X$ is the one-dimensional physical space variable,  and positive quantities $U$ and $W$ are the population densities of the prey and predator, respectively. 
 Parameter $\mathcal A>0$ is the growth factor for the prey species, $\mathcal C>0$ is the death rate for the predator without prey, $\mathcal K>0$ is the carrying capacity of the prey species. Parameters $\mathcal B>0$ and $\mathcal D>0$ are the interaction rates for the two species.; $\mathcal E_1$ and $\mathcal E_2>0$
  reflect ``satiation" effects \cite{Dunbar} for the prey and predator, or the environmental protection/adoptability to the level of predation of the predator or prey species   \cite{AADO}. Non-negative parameters $\epsilon_u$ and $\epsilon_w $ represent diffusion of prey and predator.

To non-dimensionalize  the system \eqref{e:10}, we introduce  new variables and a new set of parameters,
\begin{equation}
\begin{array}{llll}
u=\mathcal E_1U, &w=\frac{\mathcal E_1\mathcal K\mathcal B}{\mathcal A}W, & x=\sqrt{\frac{\mathcal D-\mathcal C\mathcal E_2}{\mathcal E_2}}X,&\quad t=\frac{\mathcal D-\mathcal C\mathcal E_2}{\mathcal E_2}\tau,  \\
\alpha=\frac{\mathcal C\mathcal E_1}{\mathcal D-\mathcal C\mathcal E_2}, &\gamma=\mathcal E_1\mathcal K,&\delta= \frac{\mathcal E_1 \mathcal K(\mathcal D-\mathcal C\mathcal E_2)}{\mathcal A\mathcal E_2}, 
 &\quad \eta=\frac{\mathcal E_1}{\mathcal E_2},
\end{array}
\end{equation}
and  rewrite \eqref{e:10}  as
\begin{eqnarray}\label{e:11_0}
u_{t}&=&\epsilon_uu_{xx} +\frac{1}{\delta}\left(u\left(\gamma-u\right)-\frac{uw}{1+u}\right),\notag\\
w_{t}&=&\epsilon_ww_{xx} + \frac{ w\left(u-\alpha\right)}{\eta+u}.
\end{eqnarray} 

The parameter $\eta=\frac{\mathcal E_1}{\mathcal E_2}$ is a relative protection rate and $\gamma=\mathcal E_1\mathcal K$ is  now the scaled carrying capacity of the prey. 
Obviously, when introducing  this scaling we assume that the quantity $\,\frac{\mathcal D-\mathcal C\mathcal E_2}{\mathcal E_2}=\frac{\mathcal D}{\mathcal E_2} -\mathcal C\,$ is positive. We interpret this assumption as a  reflection of the fact  that successful predation can happen only if the  measure of  interaction of the predator with the prey  relative to the environmental protection of the predator  should be higher than the death rate of the predator in absence of the prey. The higher protection rate $\mathcal E_2$ does not  necessarily always have a positive impact, as it  may reduce the effective  interaction rate of the predator and the prey and lead to the predator's ``isolation" from its food source. The quantity $\frac{\mathcal D}{\mathcal E_2}$ is   the effective  interaction rate of the predator with the prey  which reflects the  predator's population growth when it has access  to the prey. In this paper, we concentrate on the situations when $\frac{\mathcal D}{\mathcal E_2}$  is slightly higher than the death rate of the predator  in the absence of the prey, or $0<\delta \ll 1$.
\begin{eqnarray}\label{e:11_01}
u_{t}&=&\epsilon_uu_{xx} +\left(u\left(\gamma-u\right)-\frac{uw}{1+u}\right),\notag\\
w_{t}&=&\epsilon_ww_{xx} + \delta \frac{ w\left(u-\alpha\right)}{\eta+u}.
\end{eqnarray} 
The new scaling  reduces the number of parameters in the system,  but the main motivation behind it is that  the system  \eqref{e:11_0} is  amenable for the mathematical techniques that we use in the paper.   

We make an additional assumption  
$\gamma=1$, 
in order to  simplify the algebraic calculations involved in the analysis.
Although, in general,  it is a restrictive assumption,  it does not cause any loss of generality when the system is investigated for the purposes of proving  the existence of the particular set of the traveling fronts that we found.  We discuss the influence of this assumption on the results in  Section~\ref{gamma}. To summarize, 
we study the system 
\begin{eqnarray}\label{e:11}
u_{t}&=&\epsilon_uu_{xx} +\frac{1}{\delta}\left(u\left(1-u\right)-\frac{uw}{1+u}\right),\notag\\
w_{t}&=&\epsilon_ww_{xx} + \frac{ w\left(u-\alpha\right)}{\eta+u},
\end{eqnarray}
where parameters  $\alpha$ and  $\eta$ are positive and parameters $\epsilon_u$, $\epsilon_w$, and $\delta$ are nonnegative. 

The spatially homogeneous equilibria of the system of partial differential equtions \eqref{e:11} are given by the solutions $(u,w)$ of the algebraic system 
$$u(1-u)-\frac{uw}{1+u}=0,\quad \frac{w\left(u-\alpha\right)}{\eta+u}=0.$$
 The three  equilibria $(u,w)$  are
 \begin{equation}A=(\alpha,1-\alpha^2), \,\,\,\, B=(1, 0), \,\,\,\,O=(0,0). \label{equilibria}\end{equation}
 When $ \alpha\leq 1$,  all of these equilibria have nonnegative components and therefore are relevant for the population modeling.  We are interested in the case when the equilibria $A$ and $B$ are distinct, so we assume  that $ \alpha < 1$.

In this paper we construct  solutions to the system \eqref{e:10} of a special type. These are traveling fronts which are solutions that propagate  with constant velocity without changing their shape and that asymptotically connect distinct equilibria.  An important  example of an equation that supports traveling fronts is the Fisher-KPP equation. In 1937, the Fisher-KPP equation was  independently formulated  as a model  for the  geographic spread of an advantageous gene  by Fisher \cite{Fisher} and by  Kolmogorov, Petrovsky and Piskunov  \cite{KPP}. It was shown that 
  the gene frequency behaved like a  front traveling at a  speed depending on the gene's advantage. This process is captured by  the
  Fisher-KPP equation
\begin{equation} u_t =u_{xx} +au(1-u), \label{kpp_0} \end{equation}
where $a > 0$ is a parameter, $u(x,t)$ is the proportion of the population located at point $x$ at time $t$ that possesses the favorable gene. The dynamical properties of the fronts  that  asymptotically connect  $u=0$ and $u=1$ in the Fisher-KPP equation are well understood. 

A  more general  form of the Fisher-KPP equation  is
\begin{equation}u_t =u_{xx} +f(u).\label{KPP}\end{equation}
Under certain conditions on $f$,  this equation  shares the properties of the original equation \eqref{kpp_0}.  For the system under our consideration these conditions are formulated and checked  later, in \eqref{KPPnon}.
 
In this paper,  we  investigate traveling fronts  in  \eqref{e:11} in the following situations: 
\begin{itemize}
 \item[Case 1.] The quantity describing the prey diffuses at the rate much smaller than that of the predator.  
 \item[ Case 2.]   Both the predator and the prey diffuse very slowly. 
  \item[ Case 3.]  The fronts move at  fast speeds. 
   \item[ Case 4.]  The internal scaling of the system is quite large; in other words, the size of the spatially localized structure that one expects to see there is large. 
 \end{itemize}

Some of these situations are clearly related to specific modeling situations. For example, in Case  1 the assumption is natural for population models  for herbivore predation and for models describing vegetation propagation patterns. In general, all of the above situations may be described by  specific  relations between parameters which cover significant areas in the parameter space. To introduce these relations, 
in what follows, we use $\epsilon$ as a parameter that  captures smallness. Its specific definition will vary from case to case.  More precisely,  \begin{itemize}
\item In Case 1 we denote 
 $\epsilon=\frac{\epsilon_u}{\epsilon_w}\ll 1$, $\epsilon\ll\delta\ll 1$.
\item In Case 2 we consider waves that propagate with speed $c=O(1)$ in vanishing diffusion limit:  $\epsilon_u$ and  $\epsilon_w$ are of order $O(\epsilon)$. 
\item In Case 3, we
consider  fast waves with speed $c= O(\epsilon^{-1/2})$,  while $\epsilon_u$  and $\epsilon_v$ are of order $O(1)$.  
\item 
 For Case 4,  we use a scaling  of the spatial variable  that  reveals  the respective  smallness of the diffusion coefficients. 
 \end{itemize}
 
{\bf Remark}. In  the special case of  \eqref{e:10}  with   $\mathcal E_1=\mathcal E_2$, in \cite{Fu}  traveling fronts  were proved to exits  using approach from \cite{B}.  %We point out that   
 In \cite{Fu}, the authors mention that  their techniques developed  for the analysis of traveling fronts in  the Fisher-KPP scalar  equation  are not  applicable for this system  because of the lack of the comparison principle. 
 The results of the current paper actually imply that in the parameter regimes  related to  Case 1,  the dynamical properties  in the full system are strongly dominated  by  the  dynamics  of a Fisher-KPP equation \eqref{KPP}. In particular, the  cornerstone of the proof of the  existence of the traveling fronts  is  the technique used for the Fisher-KPP equation, as it will be demonstrated below.
 
We  first investigate Case 1, then we  address the existence of traveling fronts  for the parameter  regimes  described in Cases 2-4  in a unified way.

\section{Case 1. Slowly  diffusing  prey  \label{sec:1d}}

\subsection{Scaling and  formulation of the result.}
In this section we consider the situation when the prey $u$   diffuses at the rate much smaller than that of the predator $w$, i. e.,  in mathematical terms,  
$\epsilon=\frac{\epsilon_u}{\epsilon_w}\ll 1$.  We assume that  $\epsilon$ is  the smallest parameter, so  that $\epsilon\ll\delta\ll 1$.

 We  change the spatial variable in \eqref{e:11}  to
  $z = x /\sqrt{\epsilon_w} $ to obtain
\begin{eqnarray}\label{e:110}
u_{t}&=&\epsilon u_{zz} +\frac{1}{\delta}\left(u\left(1-u\right)-\frac{uw}{1+u}\right),\notag\\
w_{t}&=& w_{zz} + \frac{ w\left(u-\alpha\right)}{\eta+u},
\end{eqnarray} 
and then, to capture traveling waves, we pass to the moving frame $\zeta = z-ct$, where $c$ is a parameter representing the speed of the waves. 
When $c=0$, these stationary solutions represent standing waves. We do not study standing waves in this paper and focus only on   traveling fronts that move with with velocity $c\neq 0$.   

In the new variable $\zeta$ the equation \eqref{e:110}  reads
\begin{eqnarray}\label{e:21}
u_{t}&=&{\epsilon} u_{\zeta\zeta} +c u_{\zeta}+ \frac{1}{\delta}\left(u\left(1-u\right)-\frac{uw}{1+u}\right),\notag\\
w_{t}&=&  w_{\zeta\zeta}+ c w_{\zeta}+\frac{ w\left(u-\alpha\right)}{\eta+u}.
\end{eqnarray}   
The system \eqref{e:21} is  invariant under the transformation $(\zeta,c)\to (-\zeta, -c)$,  and, thus, it is enough  to consider $c>0$. 

The main result of this section is the following theorem. 
\begin{Theorem} \label{T:1}
For every fixed  $0<\alpha<1$, $\eta >0$,  and every  $c > 2\sqrt{\frac{1-\alpha}{\eta+1}  }$,  there exists  $\delta_0=\delta_0(\alpha,\eta,c) >0$    such that for every $0<\delta <\delta_0$  there exists $\epsilon_0(\alpha, \eta,c,\delta)>0$  such that for each $0<\epsilon <\epsilon_0$  there is a  translationally invariant family of  fronts  of the system \eqref{e:110} which   move with speed $c$,  converge to the  equilibrium  $A=(\alpha, 1-\alpha^2)$  at $-\infty$  and  to   the equilibrium   $B =(1,0)$ at $+\infty$, and, moreover,  which have positive components $u$ and $w$.  
\end{Theorem}

\subsection{ Analysis of the traveling wave system.}
Traveling fronts  that move with velocity $c$ are  stationary solutions of \eqref{e:21}, in other words, they are solutions of a system 
\begin{eqnarray}\label{e:tw}
0&=&{\epsilon} u_{\zeta\zeta} +c u_{\zeta}+ \frac{1}{\delta}\left(u\left(1-u\right)-\frac{uw}{1+u}\right),\notag\\
0&=&  w_{\zeta\zeta}+ c w_{\zeta}+\frac{ w\left(u-\alpha\right)}{\eta +u},
\end{eqnarray} 
which is  a system of ordinary differential equations. 
We  use the coordinate transformation
\begin{equation} \label{e:cord} u_1=u, \quad u_2=u_\zeta,\quad w_1=w,\quad w_2=w_\zeta, \end{equation}
 to rewrite the traveling wave system \eqref{e:21} as a system of  the first order ordinary differential equations 
\begin{eqnarray}\label{e:31_0}
\frac{du_1}{d\zeta}&=&u_2,\notag\\
\epsilon\frac{du_2}{d\zeta}&=&-c u_2+\frac{1}{\delta}\left(\frac{u_1w_1}{1+u_1}- u_1(1-u_1)\right),\notag\\
\frac{dw_1}{d\zeta}&=&w_2,\notag\\
\frac{dw_2}{d\zeta}&=&-c w_2-{}  \frac{w_1\left(u_1-\alpha\right)}{\eta+u_1}.
\end{eqnarray} 
In addition, we consider this system when the independent variable is rescaled as   $\xi =\zeta/ \epsilon$. In terms of the variable $\xi$ the system \eqref{e:31_0} reads
\begin{eqnarray}\label{e:31_0f}
\frac{du_1}{d\xi}&=&\epsilon u_2,\notag\\
\frac{du_2}{d\xi}&=&-c u_2+\frac{1}{\delta}\left(\frac{u_1w_1}{1+u_1}- u_1(1-u_1)\right),\notag\\
\frac{dw_1}{d\xi}&=& \epsilon w_2,\notag\\
\frac{dw_2}{d\xi}&=&\epsilon\left(-c w_2- \frac{w_1\left(u_1-\alpha\right)}{\eta+u_1}\right).
\end{eqnarray} 
The system \eqref{e:31_0} is called the slow system and the system \eqref{e:31_0f} is called the fast system. 

  The system  \eqref{e:31_0f} can be considered a multi-scale dynamical system. Indeed,  the scale separation between $\zeta$ and $\xi$ is caused by the smallness of $\epsilon$. In addition, we assume that $\delta$ is also a small parameter, so  there is an additional scale related to the smallness of $\delta$.  Multi-scale slow-fast systems arise very often in various applications.  %General results for treatment of some classes of multi-scale systems is obtained in \cite{CT}, where  Fenichel's theorems \cite{Fenichel79}  were extended   for  multi-scale slow-fast systems. In the system \eqref{e:31_0f} the scale separation is structured differently and the results from \cite{CT} do not directly apply. \newline 
General results for treatment of some classes of multi-scale systems exist, for example:
 Fenichel's theorems \cite{Fenichel79}  are extended  in \cite{CT}  for  multi-scale slow-fast systems  which  in the case of two parameters are of the  form $x_ 1^{\prime} = f_1(x; \epsilon_1; \epsilon_2)$, $x_ 2^{\prime} = \epsilon_1 f_2(x; \epsilon_1; \epsilon_1)$, $x_ 3^{\prime}  = \epsilon_1\epsilon_2 f_3(x; \epsilon_1;\epsilon_2)$.  In the system \eqref{e:31_0f} the scale separation is structured differently and the results from \cite{CT} do not directly apply. 
 
   We  assume that  $\epsilon\ll\delta\ll 1$,  and  use this assumption  to reduce the system \eqref{e:31_0} using the  geometric singular perturbation  theory twice, one time with respect to the small parameter $\epsilon$ and the second time with respect to $\delta$  ($\epsilon\ll\delta\ll 1$). This type of approach  was used in a number of studies of multi-scale systems. For example,  in  \cite{M_GS} for the geometric construction of traveling waves in the Gray-Scott model, in \cite{GM}  for the construction of traveling waves in a population model for the mussel-algae interaction, and fronts and periodic traveling waves in a diffusive Holling-Tanner population model \cite{GMS}. 

 %\textcolor{blue}{  The system  \eqref{e:31_0f} can be considered a multi-scale dynamical system. Indeed,  the scale separation between $\zeta$ and $\xi$ is caused by the smallness of $\epsilon$. In addition, we assume that $\delta$ is also a small parameter, so  therefore there is an additional scale related to the smallness of $\delta$.  Multi-scale slow-fast systems arise very often in various applications. 

 %We  assume that  $\epsilon\ll\delta\ll 1$,  and  use this assumption  to reduce the system \eqref{e:31_0} using the  geometric singular perturbation  theory, twice, one time with respect to the small parameter $\epsilon$ and the second time with respect to $\delta$ once   with respect to the small parameter $\epsilon$ and the second time with respect to $\delta$ ($\epsilon\ll\delta\ll 1$). The similar approach was used in a number of studies of multi-scale systems, for example,  in  \cite{M_GS} for the geometric construction of traveling waves in Gray-Scott model, in \cite{GM}  for the geometric construction of traveling waves in a population model for mussel-algae interaction, or in \cite{GMS} for construction of fronts in a diffusive Holling-Tanner population model. }
 
%  the assumption  $\epsilon\ll\delta\ll 1$  to reduce the system \eqref{e:31_0} using the  geometric singular perturbation  theory, twice, one time with respect to the small parameter $\epsilon$ and the second time with respect to $\delta$.

When we consider the limit  of  \eqref{e:31_0}  as  $\epsilon \to 0 $,   the second equation yields an algebraic relation for the components of the solution $(u_1,u_2,w_1,w_2)$, which defines the 3-dimensional  slow manifold: 
\begin{equation}M_{\epsilon=0}=\left\{(u_1,u_2,w_1,w_2): u_2=\frac{1}{c\delta}\left(\frac{u_1w_1}{1+u_1}- u_1(1-u_1)\right)\right\}.\end{equation}
On the set  $M_{\epsilon=0} $, the system \eqref{e:31_0} is reduced to 
\begin{eqnarray}\label{e:31_b}
\delta\frac{du_1}{d\zeta}&=&\frac{1}{c}\left(\frac{u_1w_1}{1+u_1}- u_1(1-u_1)\right),\notag\\
\frac{dw_1}{d\zeta}&=&w_2,\notag\\
\frac{dw_2}{d\zeta}&=&-c w_2-{}  \frac{w_1\left(u_1-\alpha\right)}{\eta+u_1}.
\end{eqnarray} 

The slow manifold   $M_{\epsilon=0} $  also represents the set of equilibria for \eqref{e:31_0f} when $\epsilon=0$
\begin{eqnarray}\label{e:31_0f0}
\frac{du_1}{d\xi}&=&0,\notag\\
\frac{du_2}{d\xi}&=&-c u_2+\frac{1}{\delta}\left(\frac{u_1w_1}{1+u_1}- u_1(1-u_1)\right),\notag\\
\frac{dw_1}{d\xi}&=& 0, \notag\\
\frac{dw_2}{d\xi}&=&0.
\end{eqnarray} 
  The linearization of the  system \eqref{e:31_0f0}  at each point of $M_{\epsilon=0} $  has three  zero eigenvalues and one  negative eigenvalue $-c$, therefore $M_{\epsilon=0} $  is   normally hyperbolic and attracting. By Fenichel's invariant manifold theory \cite{Fenichel79, Jones94}, $M_{\epsilon=0} $  persists as an invariant set in the system with sufficiently small positive  $\epsilon$ in a form of an $\epsilon$-order perturbation $M_{\epsilon}  $,
  \begin{equation}M_{\epsilon}=\left\{(u_1,u_2,w_1,w_2): u_2=\frac{1}{c\delta}\left(\frac{u_1w_1}{1+u_1}- u_1(1-u_1) +O(\epsilon)\right)\right\}, \end{equation}
   with the flow on $M_{\epsilon} $ being an $\epsilon$-order perturbation of the flow \eqref{e:31_b}, 
  \begin{eqnarray}\label{e:31_be}
\delta\frac{du_1}{d\zeta}&=&\frac{1}{c}\left(\frac{u_1w_1}{1+u_1}- u_1(1-u_1)\right) +O(\epsilon),\notag\\
\frac{dw_1}{d\zeta}&=&w_2,\notag\\
\frac{dw_2}{d\zeta}&=&-c w_2-{}  \frac{w_1\left(u_1-\alpha\right)}{\eta+u_1}.
\end{eqnarray} 

%We intend to use geometric singular perturbation theory  to extend the information about orbits that exist in the limiting system \eqref{e:31_b}  to  the system \eqref{e:31_0}, but 
 We next analyze the system \eqref{e:31_b} using the smallness of $\delta$.  We consider \eqref{e:31_b} simultaneously with its rescaled version with respect to the new variable $\varsigma=\zeta /\delta$, 
\begin{eqnarray}\label{e:31_b_delta}
\frac{du_1}{d\varsigma}&=&\frac{1}{c}\left(\frac{u_1w_1}{1+u_1}- u_1(1-u_1)\right),\notag\\
\frac{dw_1}{d\varsigma}&=&\delta w_2,\notag\\
\frac{dw_2}{d\varsigma}&=&\delta\left(-c w_2- \frac{w_1\left(u_1-\alpha\right)}{\eta+u_1}\right).
\end{eqnarray} 

In  the singular limit of \eqref{e:31_b} as $\delta\to 0$,  
\begin{eqnarray}\label{e:31_b0}
0 &=&\frac{1}{c}\left(\frac{u_1w_1}{1+u_1}- u_1(1-u_1)\right),\notag\\
\frac{dw_1}{d\zeta}&=&w_2,\notag\\
\frac{dw_2}{d\zeta}&=&-c w_2-  \frac{w_1\left(u_1-\alpha\right)}{\eta+u_1}.
\end{eqnarray} 
the solutions    live on  $M_{\delta=0}^1 \cup M_{\delta=0}^2$, where 
 \begin{equation} \label{MM}M_{\delta=0}^1=\{(u_1,w_1,w_2): w_1=1-u_1^2\}\quad
\text{ and } \quad  M_{\delta=0}^2 = \{(u_1,w_1,w_2): u_1=0\}.\end{equation}
 On $M_{\delta=0}^2$  the dynamics of the system \eqref{e:31_b0} is that of the linear equation 
$$\frac{d^2 w_1}{d\zeta^2 } +c \frac{dw_1}{d\zeta} -\frac{\alpha}{\eta} w_1=0.$$ 
 No solution of this equation converges to constant states at both $\pm \infty$, therefore it is of no relevance to the analysis.
 For $M_{\delta=0}^1$, we consider only the portion of it  where  $u_1=u >0$ %and $w_1=w \geq 0$, so we denote
 $$M_{\epsilon=0, \delta=0}= \{(u_1,w_1,w_2): u_1=\sqrt{1-w_1},  w_1 <1\}.$$
 
 The set $M_{\epsilon=0, \delta=0}$ is a  portion of the set of  all equilibria of  \eqref{e:31_b_delta} with $\delta=0$,
 \begin{eqnarray}\label{e:31_b_delta_0}
\frac{du_1}{d\varsigma}&=&\frac{1}{c}\left(\frac{u_1w_1}{1+u_1}- u_1(1-u_1)\right),\notag\\
\frac{dw_1}{d\varsigma}&=&0,\notag\\
\frac{dw_2}{d\varsigma}&=&0.
\end{eqnarray}
The linearization of  \eqref{e:31_b_delta_0} around any point of the two-dimensional  set $M_{\epsilon=0, \delta=0}$ has two zero eigenvalues and  the  eigenvalue $\frac{2(1-w_1)}{c(1+\sqrt{1-w_1})}$. This eigenvalue is strictly positive when $w_1<1$, and so the set is  normally hyperbolic and repelling.  

By Fenichel's First Theorem, the critical manifold $M_{\epsilon=0, \delta=0}$, at least over compact sets, perturbs to an invariant manifold $M_{\epsilon=0, \delta}$ for \eqref{e:31_b_delta} with $\delta > 0$ but sufficiently small, and $M_{\epsilon=0, \delta}$ is  a  $\delta$-order perturbation of  $M_{\epsilon=0, \delta=0}$, where  $u_1=\sqrt{1-w_1}+O(\delta)$.  If $\delta$ is small enough,  $M_{\epsilon=0, \delta}$ is also  normally hyperbolic and repelling on the fast scale  $\varsigma=\zeta/\delta$.

On $M_{\epsilon=0, \delta}$  the flow on the slow scale $\zeta$ is given by $\delta$-order perturbation of  the flow on $M_{\epsilon=0, \delta=0}$.
More precisely, the  flow on $M_{\epsilon=0, \delta=0}$ is given by
the reduced system:
\begin{eqnarray}\label{e:31_r}
\frac{dw_1}{d\zeta}&=&w_2,\notag\\
\frac{dw_2}{d\zeta}&=&-c w_2- \frac{w_1\left(\sqrt{1-w_1}-\alpha\right)}{\eta+\sqrt{1-w_1}}.
\end{eqnarray} 
The flow on $M_{\epsilon=0, \delta}$  is a $\delta$-order perturbation of the \eqref{e:31_r},
\begin{eqnarray}\label{e:31_rd}
\frac{dw_1}{d\zeta}&=&w_2,\notag\\
\frac{dw_2}{d\zeta}&=&-c w_2- \frac{w_1\left(\sqrt{1-w_1}-\alpha\right)}{\eta+\sqrt{1-w_1}}  +O(\delta).
\end{eqnarray} 

We are interested in  a possibility of  heteroclinic orbits of \eqref{e:31_b_delta} and  \eqref{e:31_b_delta_0} that asymptotically connect equilibria $(1,0,0) $ and $(\alpha, 1-\alpha^2, 0)$ that correspond to equilibria  $B $  and  $A$  in \eqref{equilibria}. Thus, we  concentrate on  equilibria $\tilde B= (w_1,w_2)=(0,0)$   and $\tilde A =(1-\alpha^2,0)$ for 
the system \eqref{e:31_r}. 
 
We note that  at $\alpha =1$,  the reduced system \eqref{e:31_r} undergoes a saddle-node bifurcation. For a  heteroclinic orbit  to exist it is critical to have two distinct  equilibria with non-negative components so we assume that   $0<\alpha<1$.  

%By Fenichel's Theory any invariant set  for \eqref{e:31_b_delta} which is sufficiently close  to $M_{\epsilon=0, \delta=0}$  is located on $M_{\epsilon=0, \delta}$ , thus, since t
The equilibria $(u_1,w_1,w_2)= (\alpha,1-\alpha^2,0)$ 
 and $(u_1,w_1,w_2)= (1, 0, 0)$  which in the system \eqref{e:11} correspond to the equilibria $A$ and $B$, respectively, belong to  %$M_{\epsilon=0, \delta=0}$ and 
 $M_{\epsilon=0, \delta}$ as well, since they are equilibria of   the system %each of the systems\eqref{e:31_r} and
 \eqref{e:31_b_delta}.

We show below that, as long as  $c \geq 2\sqrt{\frac{1-\alpha}{\eta+1}  }>0$ in the system \eqref{e:31_r},  the equilibrium $(0,0)$ is an attractor and  the equilibrium $(1-\alpha^2,0)$ is a saddle.  For each  such fixed $c$,  one can choose $\delta>0$ small  enough   that  these equilibria in the system \eqref{e:31_rd}  are qualitatively the same as in  \eqref{e:31_r}, an attractor and a saddle. %; to be  also an attractor and a saddle, respectively.
Moreover, we show that when $\delta=0$, the 2-dimensional stable manifold of $\tilde B =(0,0)$  intersects  the one-dimensional unstable manifold of $\tilde A=(1-\alpha^2,0)$, transversally.  
\begin{Lemma} \label{L:1}
 For every fixed $0<\alpha<1$,  $\eta >0$, and  $c \geq 2\sqrt{\frac{1-\alpha}{\eta+1}  }$, the system  \eqref{e:31_r} has a heteroclinic orbit   that asymptotically connects     the saddle $(1-\alpha^2,0)$  at $-\infty$  to  the node  $(0,0)$ at $+\infty$.
 \end{Lemma}
 \begin{Proof} 
The linearization of \eqref{e:31_r} about $(1-\alpha^2,0)$  has  two  real
 eigenvalues  $$\lambda_{1,2}=\frac{-c\pm \sqrt{c^2+2\frac{1-\alpha^2}{\alpha(\eta+\alpha)}}}{2}$$ of opposite signs, with eigenvectors $(1,\lambda_{1,2})$,   so $(1-\alpha^2,0)$  is a saddle. 

The linearization of \eqref{e:31_r} about $(0,0)$  has  two 
 eigenvalues  $$\lambda_{1,2}=\frac{-c\pm \sqrt{c^2 -4\frac{1-\alpha}{\eta+1}}}{2},$$ 
 %both of which are negative, 
 with corresponding eigenvectors $(1,\lambda_{1,2})$. %, so $(0,0)$ is a node.
 Under the assumption $0<\alpha<1$,   $(0,0)$ is stable node: the eigenvalues are real and negative as long as  $c^2-4\frac{(1-\alpha)}{\eta+1}>0$; they are 
 complex conjugate with negative real parts when   $c^2-\frac{4(1-\alpha)}{\eta+1} <0$. The critical value of $c =2\sqrt{\frac{1-\alpha}{\eta+1}  }$ captures the transition from oscillatory to monotone convergence to the equilibrium.

To argue the  existence of heteroclinic orbits that connect $(0,0)$ and $(1-\alpha^2,0)$  we point out that the system \eqref{e:31_r}  is equivalent to the scalar ode corresponding to the traveling wave equation for a Fisher-KPP  equation  \eqref{KPP}.
Indeed, \eqref{e:31_r} is equivalent to
$$\frac{d^2w_1}{d\zeta^2 }+c \frac{dw_1}{d\zeta} + \frac{w_1\left(\sqrt{1-w_1}-\alpha\right)}{\eta+\sqrt{1-w_1}} =0.$$
The function
 $$f(w_1)=  \frac{w_1\left(\sqrt{1-w_1}-\alpha\right)}{\eta+\sqrt{1-w_1}}$$ 
satisfies the conditions on Fisher-KPP nonlinearity \cite{Fisher, KPP}: 
\begin{eqnarray} 
&f(0)=f(1-\alpha^2)=0, \quad f^{\prime} (0)=\frac{1-\alpha}{\eta+1} >0,\quad  f^{\prime} (1-\alpha^2) =\frac{(\alpha^2-1)}{2\alpha(\alpha+\eta)}<0, \notag \\& f^{\prime\prime}(w_1)=\frac{((3w_1-4)\eta+\sqrt{1-w_1}(w_1-4))(\alpha+\eta)}{4(1-w_1)^{3/2}(\eta+\sqrt{1-w_1})^3}<0. \label{KPPnon}\end{eqnarray}
The latter inequality holds because  $0\leq w_1<1$ on $M_{\epsilon=0, \delta=0}$.

 \begin{figure}[t]
\begin{center}
\includegraphics[width=3.in]{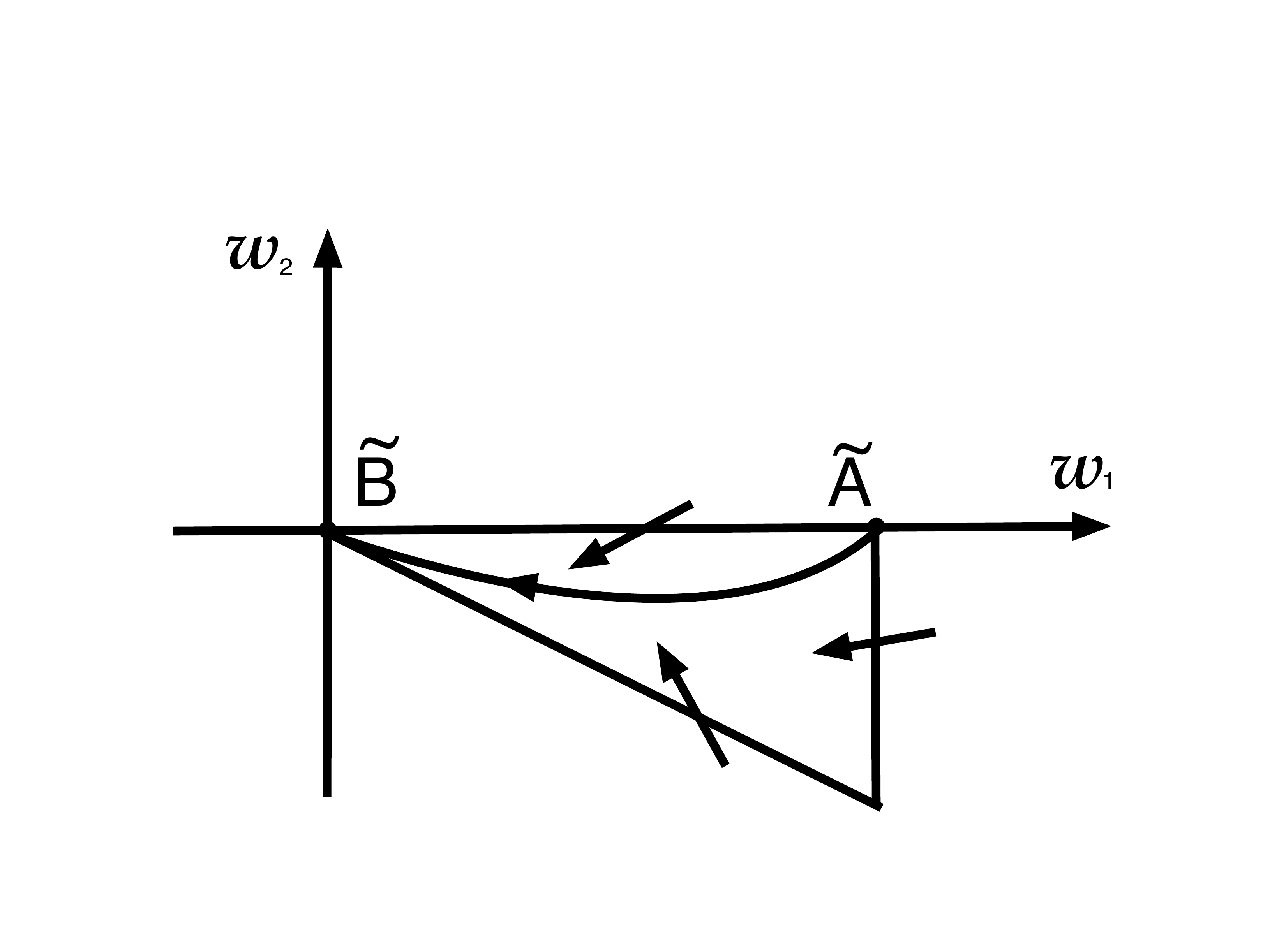} 
\caption{Existence of a heteroclinic orbit for\eqref{e:31_r}  by means of  a trapping region argument.}
\label{fig:trapping}
\end{center}
\end{figure}

 It is known \cite{KPP, AW1,AW2}  that for   $c\geq 2\sqrt {f^{\prime} (0)}=2\sqrt{ \frac{1-\alpha}{\eta+1}}$,  a monotone in $w_1$  asymptotic connection between the equilibria exists. %The monotonicity in $w_1$ implies that  $w_1>0$.
 
The classical proof  of this fact  is based on a construction of a trapping region  in  $(w_1,w_2)$-space for solutions of \eqref{e:31_r}. 
 In our case, the trapping region (see Fig.~\ref{fig:trapping}) is bounded by  a  triangle  that consists of a segment of the $w_1$-axis between $(0,0)$ and $(1-\alpha^2, 0)$, a segment of the vertical line  and a segment along line $w_2=-b w_1$, 
where $b$ is a number described below. The vector field points inside that triangular region through the vertical and horizontal sides. Along the side $w_2=-bw_1$,  for the slope of the vector field we have  %is steeper than that of the line  $w_2=-bw_1$.  Indeed, 
$$\frac{d w_2}{d w_1} = \frac{ -c w_2-   \frac{w_1(\sqrt{1-w_1}-\alpha)}{\eta+\sqrt{1-w_1}}} {w_2}\big |_{w_2=-bw_1} = \frac{ c bw_1-   \frac{w_1(\sqrt{1-w_1}-\alpha)}{\eta+\sqrt{1-w_1}}}{-bw_1}=-c +\frac{ (\sqrt{1-w_1}-\alpha)}{b(\eta+\sqrt{1-w_1})}.$$
Since the function $\frac{ s-\alpha}{\eta+s}$ is increasing, then over the interval $[ \alpha, 1]$ its maximum  value  is $\frac{ 1-\alpha}{\eta+1}$,  so
$$\frac{d w_2}{d w_1}  \leq -c+ \frac{  (1-\alpha)}{b(\eta+1)}.$$
If  $c > 2 \sqrt{\frac{1-\alpha}{\eta+1}}$ and  $b$ is a number strictly between $\frac {1}{2}\big(c\pm \sqrt{c^2-\frac{4(1-\alpha)}{\eta+1}}\big)$ then 
$c+ \frac{  (1-\alpha)}{b(\eta+1)} <-b $ and, therefore,
$\frac{d w_2}{d w_1} % \leq -c+ \frac{  (1-\alpha)}{b(\eta+1)} 
<-b$.
%%%%%%%%%%%%%%%%%
%\textcolor{blue}{$$c+ \frac{  (1-\alpha)}{b(\eta+1)} <-b \quad \text{iff}\quad  b^2-cb+ \frac{  (1-\alpha)}{b(\eta+1)}<0$$}
%%%%%%%%%%%%%%%%
Thus, the vector field points inside the triangular region. 
If  $c = 2 \sqrt{\frac{1-\alpha}{\eta+1}}$ and $b= \frac {c}{2}$, then 
$\frac{d w_2}{d w_1}  \leq %-c+ \frac{  1-\alpha}{b(\eta+1)}  =
 -b,$
so the vector field is aligned with the side  $w_2=-bw_1$ and therefore the solution that starts inside of the triangle cannot leave  either.  
Since  the unstable manifold of $(1-\alpha^2,0)$ points into this triangular region and the orbits are monotone in $w_1$ component, then the trajectory that follows the unstable manifold  asymptotically connects to $(0,0)$. 
\end{Proof}

%Note that for the asymptotic  connections from $(0,0)$ to $(1-\alpha^2,0)$ the condition $ 1-w_1>0$ is satisfied.

Next we show that this limiting orbit persists  for sufficiently small $\delta>0$, and,  then we argue that for every such $\delta$,  there exists a sufficiently small $\epsilon >0$ such that the orbit persists  in the full system \eqref{e:31_0}. 

By the dimension counting, in the two-dimensional phase space, the one-dimensional unstable manifold   of $(1-\alpha^2,0)$  and the two-dimensional stable manifold   of $(0,0)$  intersect  transversally. 
%Now we go back to the system  \eqref{e:31_b_delta}  with  sufficiently small $\delta>0$.
 In the three-dimensional  dynamical system \eqref{e:31_b_delta} with  sufficiently small $\delta>0$  the equilibrium $(1,0,0)$ has a two-dimensional stable manifold and a one-dimensional unstable manifold.   Since  $M_{\epsilon=0, \delta}$ is repelling, any solution  of \eqref{e:31_b_delta} approaching  $(1,0,0)$  does so along $M_{\epsilon=0, \delta}$ and must  entirely belong to it.    On the other hand,    the equilibrium $(\alpha,1-\alpha^2,0)$  has a two-dimensional  unstable manifold  and a one-dimensional stable manifold. 
 Moreover, this  two-dimensional  unstable manifold  includes a direction transversal to  $M_{\epsilon=0, \delta}$.  Therefore the intersection of the  two-dimensional  unstable manifold of $(\alpha,1-\alpha^2,0)$  with $M_{\epsilon=0, \delta}$ and thus with the  two-dimensional stable manifold of  $(1,0,0)$  is transversal.  This intersection forms a heteroclinic orbit for 
 \eqref{e:31_b_delta} with  sufficiently small $\delta>0$.   %The   heteroclinic orbit which is a small perturbation of the limiting orbit  with $\delta=0$  belongs to $M_{\epsilon=0, \delta}$ as an intersection of the stable manifold of   In the limit $\delta\to 0$, the limiting sets  of these two invariant manifolds  intersect transversally.
% By dimension counting ($2$ and $2$ in $3$-dimensional space) the intersection is transversal. Therefore  it persists for sufficiently small $\delta >0$. 
 We have the following lemma.

\begin{Lemma} \label{L:22} For every fixed  $0<\alpha<1$, $\eta >0$,  and $c > 2\sqrt{\frac{1-\alpha}{\eta+1}  }$, 
there is $\delta_0 =\delta_0(\alpha,\eta,c)>0$ such that  for any $0<\delta<\delta_0$  there exists a  heteroclinic orbit  of system \eqref{e:31_b} (equivalently,  \eqref{e:31_b_delta})
that  asymptotically connects     the saddle  $(u_1,w_1,w_2)=(\alpha, 1-\alpha^2,0)$  at $-\infty$  to  the saddle  $(1,0,0)$ at $+\infty$.
\end{Lemma}

Next we  consider \eqref{e:31_0f} or equivalently  \eqref{e:31_0}.  In the limit $\epsilon\to 0$, the analysis of  the flow on $M_{\epsilon=0}$  resulted in Lemma~\ref{L:22}.  
The equilibria $(\alpha, 1-\alpha^2,0)$  and   $(1,0,0)$  from \eqref{e:31_b_delta} lie on   $M_{\epsilon=0}$ and  they  are also equilibria  of  the perturbed flow  generated by  \eqref{e:31_be} on the set $M_{\epsilon}$.

  For sufficiently small $\epsilon$ the nature of equilibria  is the same, one of them is still a saddle  and the other is still a stable node.  Since $M_{\epsilon}$ is attracting,  the unstable manifold of the saddle $ (\alpha,0, 1-\alpha^2,0)$  in  \eqref{e:31_0f}  or \eqref{e:31_0} stays on $M_{\epsilon}$ as $\zeta$ changes from $-\infty$ to $+\infty$. Therefore the orbit is formed by the intersection of  the unstable manifold of the saddle $ (\alpha,0, 1-\alpha^2,0)$ and the stable manifold  of the saddle $(1,0,0,0)$.  This intersection persists since the intersection of the limits  of involved invariant sets within $M_{\epsilon=0}$ is transversal and and  entirely belongs to $M_{\epsilon}$. 
This argument proves the following theorem.
\begin{Theorem} \label{T:sdp}
For every fixed  $0<\alpha<1$, $\eta >0$,  and $c > 2\sqrt{\frac{1-\alpha}{\eta+1}  }$,  there exists  $\delta_0=\delta_0(\alpha,\eta,c) >0$    such that for every $0<\delta <\delta_0$  there exists $\epsilon_0(\alpha, \eta,c,\delta)>0$  such that for each $0<\epsilon <\epsilon_0$  there is a  heteroclinic orbit  of system \eqref{e:31_0} (equivalently,  \eqref{e:31_0f})
that  asymptotically connects     the saddle  $(u_1,u_2,w_1,w_2)=(\alpha, 0,1-\alpha^2,0)$  at $-\infty$  to  the saddle  $(u_1,u_2,w_1,w_2)=(1,0,0,0)$ at $+\infty$.
\end{Theorem}

The dynamical system representation of a traveling front  solution of a partial differential equation  is a heteroclinic orbit, therefore for the system of partial differential equations \eqref{e:110} this theorem implies Theorem~\ref{T:1}.

{\bf Remark.} Both $u_1$ and $w_1$ components of the constructed limiting  heteroclinic orbit  of  \eqref{e:31_r}  are positive.  It follows from the construction, that for sufficiently small $\delta$ and $\epsilon$,  both $u_1$ and $w_1$ components of  heteroclinic orbit of the system  \eqref{e:31_0}  are also  positive.

%%%%%%%%%%%%
\section{Vanishing diffusion limit \label{sec:2d}}
%%%%%%%%%%%%%%
%%%%%%%%%%%%

\subsection{Scaling and parameter regimes}
In this section,  we introduce  a different  scaling  for the system  \eqref{e:11} to analyze Cases 2, 3 and 4. 

 To describe Cases 2  and 3, we set $t=\delta \tau$ and $x= \sqrt{\delta}y$ in \eqref{e:11} to obtain
\begin{eqnarray}\label{e:11_1}
u_{\tau}&=&\epsilon_u  u_{yy} +u\left(1-u\right)-\frac{uw}{1+u},\notag\\
w_{\tau}&=& \epsilon_w w_{yy} + \delta \frac{ w\left(u-\alpha\right)}{\eta+u}.
\end{eqnarray} 
In the co-moving coordinate frame  $z=y -c\tau$, the system \eqref{e:11_1} reads
\begin{eqnarray}\label{e:11_12}
u_{\tau}&=&\epsilon_u  u_{zz} + c  u_{z} +u\left(1-u\right)-\frac{uw}{1+u},\notag\\
w_{\tau}&=& \epsilon_w w_{zz} + c  w_{z}+ \delta \frac{ w\left(u-\alpha\right)}{\eta+u}.
\end{eqnarray} 
Furthermore,  we set $\zeta=z/c$,   which effectively leads to  a reduction in  the number of parameters,
\begin{eqnarray}\label{e:11_13}
u_{\tau}&=&\frac{\epsilon_u}{c^2} u_{\zeta\zeta} + u_{\zeta}+u\left(1-u\right)-\frac{uw}{1+u},\notag\\
w_{\tau}&=&\frac{\epsilon_w}{c^2}  \,w_{\zeta\zeta}+  w_{\zeta}+\frac{\delta w\left(u-\alpha\right)}{\eta+u}.
\end{eqnarray}
We assume that $\epsilon_u$ and $\epsilon_w$ are positive and set
  $\eps={\epsilon_u}/{c^2}$  and $ \mu={\epsilon_w}/{\epsilon_u}$. Clearly,  $\mu>0$  and from the assumption $c>0$ it follows that $\eps>0$. The system \eqref{e:11_13} becomes
\begin{eqnarray}\label{e:11_14}
u_{\tau}&=& \epsilon u_{\zeta\zeta} + u_{\zeta}+ u\left(1-u\right)-\frac{uw}{1+u},\notag\\
w_{\tau}&=& \epsilon  \mu \,w_{\zeta\zeta}+  w_{\zeta}+\frac{\delta w\left(u-\alpha\right)}{\eta+u}.
\end{eqnarray} 
The corresponding  traveling wave ode system then reads 
%\begin{eqnarray}\label{e:11_2}
%0&=&\epsilon_u  u_{z_1z_1} + cu_{z_1} + u\left(1-u\right)-\frac{uw}{1+u},\notag\\
%0&=& \epsilon_v w_{z_1z_1} + cv_{z_1}+  \frac{ \delta w\left(u-\alpha\right)}{\eta+u},
%\end{eqnarray} 
%Furthermore,  we set $z=z_1/c$,   which effectively leads to  a reduction in  the number of parameters,
%\begin{eqnarray}\label{e:22_2}
%0&=&\frac{\epsilon_u}{c^2} u_{zz} + u_{z}+u\left(1-u\right)-\frac{uw}{1+u},\notag\\
%0&=&\frac{\epsilon_w}{c^2}  \,w_{zz}+  w_{z}+\frac{\delta w\left(u-\alpha\right)}{\eta+u}.
%\end{eqnarray} 
%We assume that $\epsilon_u$ and $\epsilon_w$ are positive and set
  %$\eps={\epsilon_u}/{c^2}$  and $ \mu={\epsilon_w}/{\epsilon_u}$. Clearly,  $\mu>0$  and from the assumption $c>0$ it follows that $\eps>0$. The system \eqref{e:22_2} becomes
\begin{eqnarray}\label{e:24}
0&=& \epsilon u_{\zeta\zeta} + u_{\zeta}+ u\left(1-u\right)-\frac{uw}{1+u},\notag\\
0&=& \epsilon  \mu \,w_{\zeta\zeta}+  w_{\zeta}+\frac{\delta w\left(u-\alpha\right)}{\eta+u}.
\end{eqnarray} 
We assume that  $0<\eps\ll 1$ or,  in other words, that  the diffusion coefficient $\epsilon_u $ is small relative to the wave velocity $c$. This assumption is satisfied for any fixed $c$  when   $\epsilon_u$ is sufficiently small  or for   for any fixed 
$\epsilon_u$ when  $c$ is sufficiently large, thus the  system \eqref{e:24}  represents the traveling wave system for both Cases 2 and 3. 

On the other hand,  to describe  Case 4, it is enough to rescale 
 $t=\delta \tau$ and $x= \sqrt{\delta} \sqrt {d} \,y$ in \eqref{e:11}, 
 assuming that  the characteristic length $d>0$  of the spatial variable  $ \sqrt{\delta} y$ is  very large \cite{M}. Indeed,  then 
\begin{eqnarray}\label{e:11_2}
u_{\tau}&=&\frac{\epsilon_u}{d}  u_{ yy} +u\left(1-u\right)-\frac{uw}{1+u},\notag\\
w_{\tau}&=&\frac{ \epsilon_v}{d}  w_{yy} + \delta \frac{ w\left(u-\alpha\right)}{\eta+u}.
\end{eqnarray} 
In the co-moving frame $z=y-ct$, \eqref{e:11_2} is
\begin{eqnarray}\label{e:11_23}
u_{\tau}&=&\frac{\epsilon_u}{d}  u_{ zz}+cu_{z} +u\left(1-u\right)-\frac{uw}{1+u},\notag\\
w_{\tau}&=&\frac{ \epsilon_v}{d}  w_{zz} +cw_{z}+ \delta \frac{ w\left(u-\alpha\right)}{\eta+u}.
\end{eqnarray}
If we further  set  $\zeta=z/c$ and 
$\eps={\epsilon_u}/{ d c^2}$  and $ \mu=\epsilon_w/\epsilon_u$,  we obtain a system that looks exactly like \eqref{e:11_14}.
 In the rest of this section we prove the following result. 
\begin{Theorem} \label{T:2}
For every fixed  $0<\alpha<1$ and  $\eta >0$ and $c\neq 0$,  there exists  $\delta_0=\delta_0(\alpha,\eta,c) >0$    such that for every $\delta\in (0,\delta_0)$  there exists $\epsilon_0(\alpha, \eta,\delta)>0$  such that for each $0<\epsilon <\epsilon_0$  there is a  translationally invariant family of  fronts  of the system \eqref{e:11_14},  moving with speed $c$, 
which  asymptotically connect     the  equilibrium  $A=(\alpha, 1-\alpha^2)$  at $-\infty$  to   the equilibrium   $B =(1,0)$ at $+\infty$.
\end{Theorem}
We assume $c>0$ without loss of generality,  since  if a traveling wave is found for $c>0$ then a coordinate change $\tau\to -\tau$ and $\zeta\to -\zeta$ captures traveling waves with negative velocity.

\subsection{Traveling wave analysis.}

In the coordinates
$u_1=u$, $u_2=u_\zeta$, $w_1=w$, $w_2=w_\zeta$,  \eqref{e:24}  can be written as the following system of  first order  ordinary differential equations:
\begin{eqnarray}\label{e:31}
\frac{du_1}{d\zeta}&=&u_2,\notag\\
\eps\frac{du_2}{d\zeta}&=&- u_2- u_1(1-u_1)+\frac{u_1w_1}{1+u_1},\notag\\
\frac{dw_1}{d\zeta}&=&w_2,\notag\\
\eps \frac{dw_2}{d\zeta}&=&-\frac{1}{\mu} w_2-\frac{\delta}{\mu}  \frac{w_1\left(u_1-\alpha\right)}{\eta+u_1}.
\end{eqnarray}

We use geometric singular perturbation theory \cite{Fenichel79, Jones94, kuehn} to analyze \eqref{e:31}. We consider  \eqref{e:31}  together with the following   system  obtained   from \eqref{e:31} by rescaling   the independent variable as $\zeta=\eps \xi$,
\begin{eqnarray}\label{e:32}
\frac{du_1}{d\xi}&=&\eps u_2,\notag\\
\frac{du_2}{d\xi}&=&- u_2- u_1(1-u_1)+\frac{u_1w_1}{1+u_1},\notag\\\
\frac{dw_1}{d\xi}&=&\eps w_2,\notag\\
\frac{dw_2}{d\xi} &=&-\frac{1}{\mu} w_2-\frac{\delta}{\mu}  \frac{w_1\left(u_1-\alpha\right)}{\eta+u_1}.
\end{eqnarray}  
 The geometric singular perturbation theory allows  us  to extend the information obtained from  the systems \eqref{e:31} and \eqref{e:32}  when $\eps=0$  to the case of nonzero but sufficiently small $\epsilon>0$. 
 We  set  $\epsilon=0$  in \eqref{e:31}, thus obtaining an algebraic description  of a set  $M_{\epsilon=0}$ 
   to which solutions  of the limiting system   belong 
\begin{equation}\label{m0}
M_{\epsilon=0}=\left\{ (u_1,u_2, w_1, w_2):\, u_1\geq 0,\,\,  u_2=- u_1(1-u_1)+\frac{u_1w_1}{1+u_1}, \,\,w_1\geq0,\,\, w_2 =-\delta \frac{w_1\left(u_1-\alpha\right)}{\eta+u_1}\right\}
\end{equation}
 and the  system of equation   defined on $M_{\epsilon=0}$ which the solutions  of the limiting system must satisfy  
\begin{eqnarray}\label{e:31m0}
\frac{du_1}{d\zeta}&=&%f_{1}(u_1,w_1)=
\frac{u_1w_1}{1+u_1}-u_1(1-u_1),\notag\\
\frac{dw_1}{d\zeta}&=&%\delta f_{2}(u_1,w_1)=
\delta \frac{w_1\left(\alpha -u_1\right)}{\eta+u_1}.
\end{eqnarray} 

The set $M_{\epsilon=0}$  can be also described as  the set of equilibrium points  for  the limiting system  \eqref{e:32} with $\eps=0$,
\begin{eqnarray}\label{e:320}
\frac{du_1}{d\xi}&=&0,\notag\\
\frac{du_2}{d\xi}&=&- u_2- u_1(1-u_1)+\frac{u_1w_1}{1+u_1},\notag\\\
\frac{dw_1}{d\xi}&=&0,\notag\\
\frac{dw_2}{d\xi} &=&-\frac{1}{\mu} w_2-\frac{\delta}{\mu}  \frac{w_1\left(u_1-\alpha\right)}{\eta+u_1}.
\end{eqnarray} 
It is easy to see that  the two-dimensional  set  $M_{\epsilon=0}$ is normally hyperbolic and attracting.  Indeed,  besides the two zero eigenvalues, each  point of $M_{\epsilon=0}$ has  eigenvalues $-1$ and $-\frac{1}{\mu}$.  Fenichel's invariant manifold theory \cite{Fenichel79, Jones94} implies  that   $M_{\epsilon=0}$  perturbs to  an  invariant set $M_{\eps}$ for the full system,  if $\epsilon$  is sufficiently small.  $M_{\eps}$ is an   $\eps$-order perturbation of $M_{\epsilon=0}$  and   the flow on $M_\eps$  is 
 \begin{eqnarray}\label{e:31meps}
  \frac{du_1}{d\zeta}&=&\frac{u_1w_1}{1+u_1}-u_1(1-u_1) +O(\eps),\notag\\
\frac{dw_1}{d\zeta}&=&\delta \frac{w_1\left(\alpha -u_1\right)}{\eta+u_1}+O(\eps).
\end{eqnarray} 

Let $$f_{1}=\frac{u_1w_1}{1+u_1}-u_1(1-u_1), \quad f_{2}=\frac{w_1\left(\alpha -u_1\right)}{\eta+u_1}.$$
The eigenvalues of the linearization of \eqref{e:31m0} at an equilibrium are given by 
\begin{equation} \label{e:41}
\lambda_{1,2}=\frac{1}{2}\left(\delta f_{2w_1}+f_{1u_1}\pm\sqrt{(\delta f_{2w_1}+f_{1u_1})^2-4\delta(f_{1u_1}f_{2w_1}-f_{1w_1}f_{2u_1})}\right),
\end{equation}
where $f_{1}$ and $f_{2}$ and their derivatives %, $f_{iu}=\partial f_i/\partial u_1$, $f_{iw}=\partial f_i/\partial w_1$, 
 are evaluated  at the equilibrium.    The derivatives of  functions $f_{1}$ and $f_{2}$   are 
\begin{equation}\label{e:42}
\quad f_{1u_1}= \frac{w_1}{(1+u_1)^{2}}-1 + 2u_1, \quad f_{1w_1}=\frac{u_1}{u_1+1}, \quad 
f_{2u_1}= -\frac{w_1(\eta+\alpha)}{(\eta+u_1)^{2}},\quad  f_{2w_1}=\frac{\alpha -u_1}{\eta+u_1}. 
\end{equation}  
At the equilibrium $A$  
\begin{equation}\label{e:43}
f_{1u_1}(A)= \frac{2\alpha^2}{1+\alpha},\quad f_{1w_1}(A)= \frac{\alpha}{1+\alpha}, \quad
f_{2u_1}(A)= -\frac{1-\alpha^2}{\eta+\alpha},\quad f_{2w_1}(A)=0. 
\end{equation}  
Since $f_{2w_1}(A)=0$,  the eigenvalues \eqref{e:41} of the linearization of \eqref{e:31m0} at  the equilibrium $A$ are
\begin{equation} \label{e:44}
\lambda_{1,2}(A)=\frac{1}{2}\left( f_{1u_1}(A)\pm \sqrt{f_{1u_1}^2(A)+4\delta f_{1w_1}(A)f_{2u_1}(A)}\right).
\end{equation}
Since $f_{1u_1}(A)>0$, the eigenvalues $\lambda_{1,2}(A)$ are positive when $\delta\in(0,\delta_0)$, where 
\begin{equation}
\delta_0 =-\frac{f_{1u_1}^2(A)}{4 f_{1w_1}(A)f_{2u_1}(A)}
=\frac{\alpha^3(\eta+\alpha)}{(1-\alpha)(1+\alpha)^2}.
\end{equation}
Therefore, when $\delta\in (0,\delta_{0})$, the equilibrium $A$ is an unstable node.

The eigenvalues of the  linearization at the equilibrium $B=(1,0)$ are 
\begin{equation} \label{e:45}
\lambda_{1}(B)=1,\quad 
\lambda_2(B)=-\delta\frac{1-\alpha}{\eta+\alpha},
\end{equation}
therefore, when $0<\alpha<1$,  the equilibrium  $B$ is a saddle.

 \subsection{Singular front at $\delta= 0$\label{s:delta0}}
 
We are interested in proving the existence of heteroclinic  orbits of  the system \eqref{e:31m0}.  We do that by exploiting the slow-fast structure of  \eqref{e:31m0} with regards to      $\delta$ when  $\delta$  is small. 

Rescaling $\varsigma=\zeta/\delta$, we obtain the fast system
\begin{eqnarray}\label{e:31m0_fast}
\delta \frac{du_1}{d\varsigma}&=&%f_{1}(u_1,w_1)=
\frac{u_1w_1}{1+u_1}-u_1(1-u_1),\notag\\
\frac{dw_1}{d\varsigma}&=&%\delta f_{2}(u_1,w_1)=
\frac{w_1\left(\alpha -u_1\right)}{\eta+u_1}.
\end{eqnarray} 
When $\delta=0$,  the system  \eqref{e:31m0} reads
\begin{eqnarray}\label{e:4m}
\frac{du_1}{d\zeta}&=&\frac{u_1w_1}{1+u_1}-u_1(1 -u_1),\notag\\
\frac{dw_1}{d\zeta}&=&0,
\end{eqnarray} 
and the system \eqref{e:31m0_fast} reads
\begin{eqnarray}\label{e:31m0_fast_0}
0&=&\frac{u_1w_1}{1+u_1}-u_1(1-u_1),\notag\\
\frac{dw_1}{d\varsigma}&=&\frac{w_1\left(\alpha -u_1\right)}{\eta+u_1}.
\end{eqnarray} 
The system \eqref{e:4m} has the set of equilibria  which  is also  the critical manifold for \eqref{e:31m0_fast_0}  described in \eqref{MM}, but, again, only one curve out of this set  is relevant
$$M_{\epsilon=0,\,\delta =0}=\{(u_1,w_1) : u_1=\sqrt{1-w_1}, \, w_1< 1 \},$$
%which is also the critical manifold for \eqref{e:31m0_fast_0}.
 The set $M_{\epsilon=0,\,\delta =0} $ is a part of the right branch of the parabola $w_1=1-u_1^2$. The branch of the parabola with negative   $u_1$, as well as the set with $u_1=0$ do not  contain the equilibria that we are interested in.  
 
\begin{figure}[t]
\begin{center}$
\begin{array}{lclc}
\includegraphics[width=3.3in]{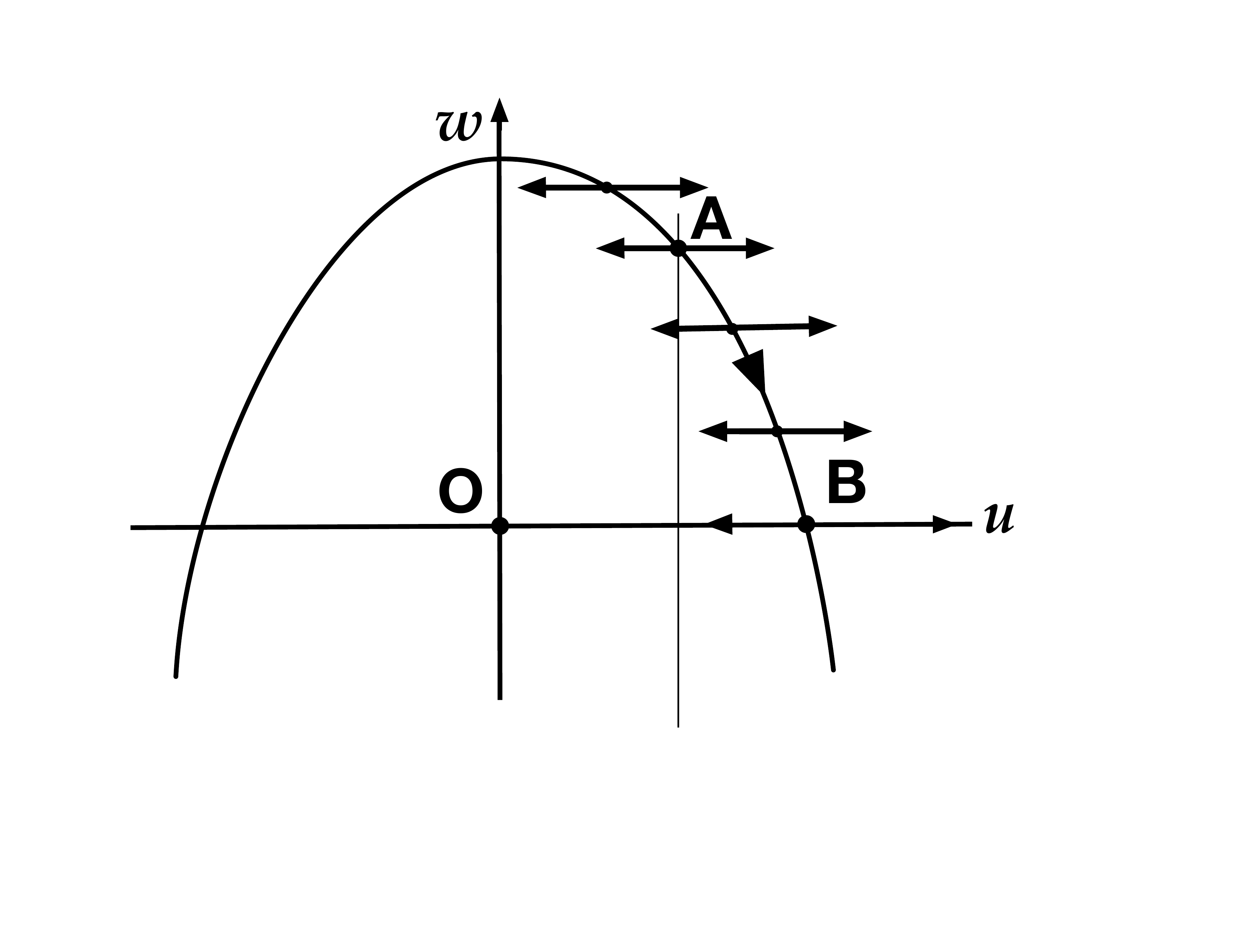} &\includegraphics[width=3.3in]{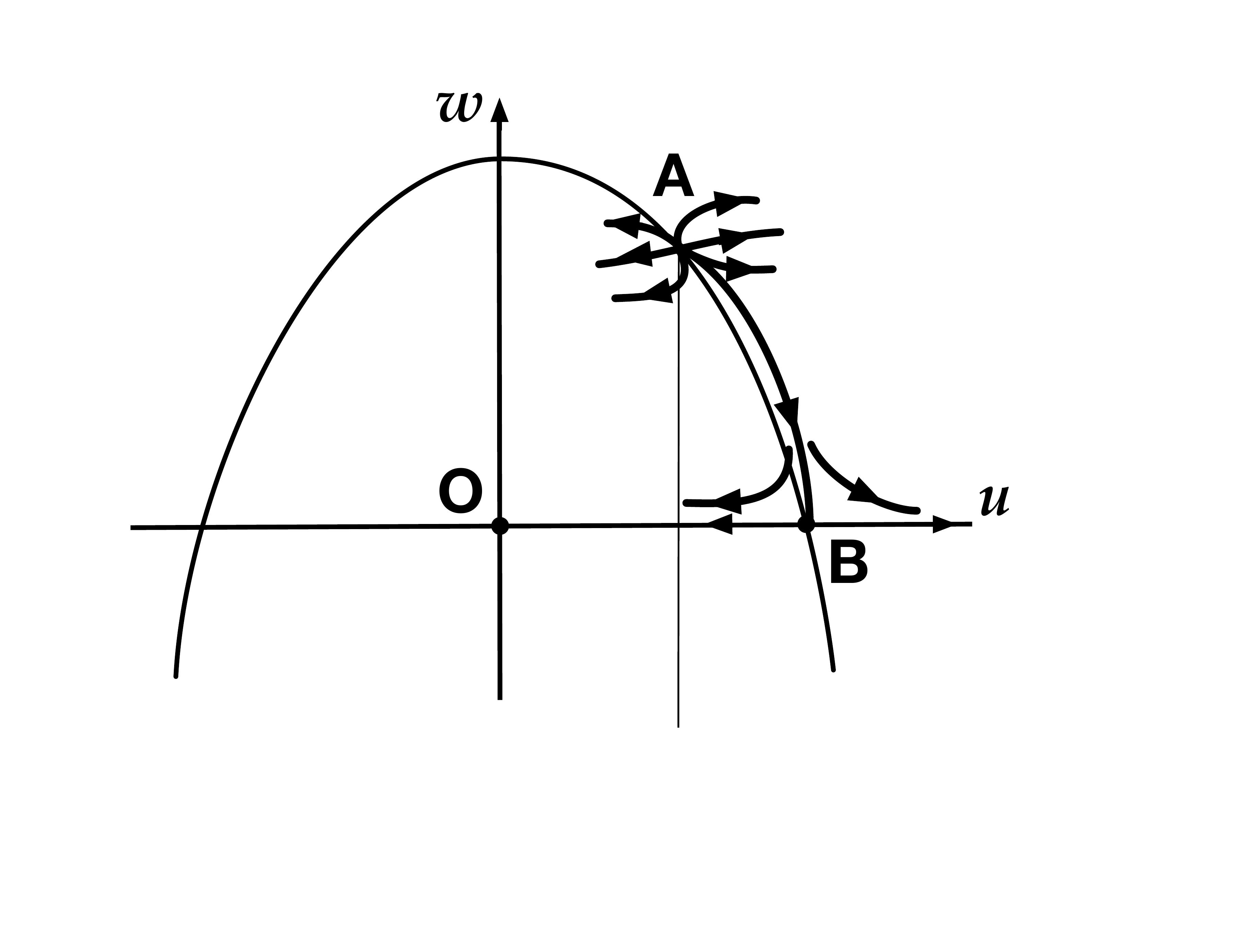} 
\end{array}$
\caption{The dynamics of \eqref{e:31m0}. Left panel: singular solution, $\delta = 0$. Right panel: perturbed solution, $0<\delta\ll 1$.}
\label{fig:waves2}
\end{center}
\end{figure}

 The linearization of \eqref{e:31m0} around $M_{\epsilon=0,\,\delta =0}$ has an eigenvalue $\frac{2(1-w_1)}{1-\sqrt{1-w_1}}$, therefore $ M_{\epsilon=0,\,\delta =0}$ is  a normally hyperbolic, repelling manifold with respect to the flow of  \eqref{e:31m0}.
On  $M_{\epsilon=0,\,\delta =0}$ the solutions  the system \eqref{e:31m0_fast_0} satisfy 
\begin{equation}\frac{dw_1}{d\varsigma}=w_1\frac{\alpha -\sqrt{1-w_1}}{\eta+\sqrt{1-w_1}},\label{slow}\end{equation}
the dynamics of which is easy to understand. 
Indeed, the  two equilibria of  the  equation \eqref{slow}  are  $w_1=0$ and $w_1=1-\alpha^2$:
  $w_1=0$  (which corresponds to  $B$)  is  stable  and  $w_1=1-\alpha^2$  (which corresponds to  $A$)  is   unstable (Fig.~\ref{fig:waves2}).

Since $M_{\epsilon=0,\,\delta =0}$ is normally hyperbolic,  it persists as an invariant set  $M_{\epsilon=0,\,\delta}$  for \eqref{e:31m0} when $\delta$ is sufficiently small.  On $M_{\epsilon=0,\,\delta}$, the flow  is a $\delta$-order perturbation of the flow described in \eqref{slow}.    Within the two-dimensional system \eqref{e:4m},   the equilibrium $A$  has  a two dimensional unstable manifold and  the equilibrium $B$ has a one dimensional stable manifold. By dimension counting their intersection is transversal, so it persists when a sufficiently small $\delta>0$ is introduced, therefore the following lemma holds.  The persistence of the  limiting orbit  that corresponds to $\delta=0$ also  follows from the general  result   \cite[Corollary 3.3]{GSz} which encapsulates   the  ``transversality  through a dimension counting" argument.  Therefore the following lemma holds.

\begin{Lemma} \label{L:31} For every fixed  $0<\alpha<1$ and  $\eta >0$,  
there is $\delta_0 =\delta_0(\alpha,\eta)>0$ such that  for any $0<\delta<\delta_0$  there exists a  heteroclinic orbit  of system \eqref{e:4m}
that  asymptotically connects  the saddle  $(u_1,w_1,w_2)=(\alpha, 1-\alpha^2,0)$  at $-\infty$  to  the saddle  $(1,0,0)$ at $+\infty$.
\end{Lemma}

This lemma proves that the singular limit of  \eqref{e:32} as $\epsilon\to 0$ has a heteroclinic connection    on  the slow manifold   $M_{\epsilon=0}$ (see \eqref{m0}).

Since   $M_{\epsilon=0}$ is normally hyperbolic and attracting it persists as an invariant manifold $M_{\epsilon}$ for the system \eqref{e:31} (equivalently, \eqref{e:32}). 

For sufficiently small $\epsilon>0$,   $M_{\epsilon}$ is attracting, therefore   the unstable manifold of the saddle $ (1,0, 0,0)$ stays on $M_{\epsilon}$ for all $\xi$. Therefore  any orbit that  originates at  $ (1,0, 0,0)$    lays on a two-dimensional set $M_{\epsilon}$, so the unstable manifold of   $ (1,0, 0,0)$ and the stable manifold of  $ (\alpha,0, 1-\alpha^2,0)$ still intersect as their limits  as $\epsilon \to 0$ intersect transversally  within $M_{\epsilon=0}$. 
The following theorem then holds. 
\begin{Theorem} \label{T:32}
For every fixed  $0<\alpha<1$ and $\eta >0$ there exists  $\delta_0=\delta_0(\alpha,\eta) >0$    such that for every $0<\delta <\delta_0$  there exists $\epsilon_0(\alpha, \eta,\delta)>0$  such that for each $0<\epsilon <\epsilon_0$  there is a  heteroclinic orbit  of the system \eqref{e:31}, or, equivalently,  of  the system \eqref{e:32})
that  asymptotically connects     the saddle  $(u_1,u_2,w_1,w_2)=(\alpha, 0,1-\alpha^2,0)$  at $-\infty$  to  the saddle  $(u_1,u_2,w_1,w_2)=(1,0,0,0)$ at $+\infty$.
\end{Theorem}

\section {The influence of the assumption $\gamma=1$. \label{gamma}}

The spatially homogeneous equilibria of the PDE \eqref{e:11_0} are
 \begin{equation}A=(\alpha,(\gamma-\alpha)(1+\alpha)), \,\,\,\, B=(\gamma, 0), \,\,\,\,O=(0,0). \label{equilibria1}\end{equation}
 When $\gamma \geq \alpha$,  all of these equilibria have nonnegative components and therefore are relevant for  population modeling.

In the limiting systems   \eqref{e:31_r} and  \eqref{e:31m0}, the co-existence equilibrium  with positive  components   corresponds to the equilibrium $A$. It is given  by the intersection of the parabolic nullcline
  $w=(1+u)(\gamma-u)$ and the verical nullcline $u=\alpha$.  We observe that the dynamics near the equilibrium depends on where the equilibrium is located relative to the vertex $V=(\frac{\gamma-1}{2},\frac{(\gamma+1)^{2}}{4})$  of the parabolic nullcline.  When the parameters are  in the region $\{ (\alpha,\gamma) : 0<\gamma< 1, \,\,\,0<\alpha <\gamma\} $   the equilibrium $A$  is in the  open first quadrant of $uw$-plane and is  strictly to the right of  the vertex  $V$.  In this situation, the proof and the  results described in  Theorems ~\ref{T:1} and \ref{T:2} hold.
The assumption $\gamma =1$  means that the vertex $V=(0,1)$ is on the vertical axis and $0<\alpha<1 $ guarantees that the equilibrium $A $ which is now $(\alpha,1-\alpha^2)$   to  the right of the   parabola's vertex.
In the parameter region $ \{(\alpha,\gamma) :  \gamma\geq 1,\,\,\,  \frac{\gamma -1}{2}<\alpha<\gamma \}$ the results still hold, although some minor, technical  changes  in the proofs are required.

\section{Holling-Tanner model \label{HT}}

The importance of the Fisher-KPP dynamics in a system has been studied before.  In  paper
\cite{Ducrot},  Ducrot considers  a diffusive predator-prey model 
\begin{eqnarray}\label{e:1h}
U_{\tau}&=&\epsilon_uU_{XX} +AU(1-\frac{U}{\mathcal K})-\frac{\mathcal BUW}{\mathcal E+U},\notag\\
W_{\tau}&=&\epsilon_w W_{XX}+\mathcal RW \left(1-\frac{\mathcal H W}{U}\right),
\end{eqnarray} 
where $t>0$  is the time and $x$ is the spacial variable,  $U$  represents the density of the prey, and $V$ is the density of the predator.  The parameters $\alpha$, $\delta$,  $\beta$, and  the diffusion coefficients $\epsilon_u$ and $\epsilon_w$ are positive.
The system \eqref{e:1h} is obtained by adding diffusion terms to an ode system for  Holling-Tanner  predator-prey model  \cite{AFPH, hsu, May0,May,M,Renshaw, T,  VP} where the predation rate is  known as Holling  type II functional response $\frac{U}{\mathcal E+U}$   \cite{H59}. 

The model \eqref{e:1h} has been extensively studied.  For example, transition fronts have been investigated numerically in \cite{UVT},  the existence of fronts in vanishing diffusion limit has been proved in \cite{GMS}, the existence of fronts for more general type of  nonlinearities was proved in \cite{A} in a parameter regime that does not capture the situation considered in this paper.  A result relevant  to our work on Rosenzweig-MacArthur system  is   the result from  \cite{Ducrot} where, in the most general case of a multi-dimensional spacial variable $x \in \mathbb R^n$, Ducrot investigates   the spreading speed of the small perturbations to  the co-existence equilibrium  which is the equilibrium with positive values for both  predator and prey densities.  The author shows in \cite{Ducrot} that  if initially  the predator is introduced in a   compactly supported manner, while  the prey is initially uniformly well distributed,  then  the spreading speed of the perturbations is defined by  the spreading speed  in the associated  Fisher-KPP scalar  equation. 

In this section we consider   \eqref{e:1h} with $x \in \mathbb R$  in a specific parameter regime and  construct fronts which are small perturbations of the fronts in scalar Fisher-KPP equation in a way  similar to one used in Section~\ref{sec:1d}. We spare readers the repetition of the details and present only  a brief overview of the analysis.

We consider a new  scaling of the system \eqref{e:1h}, given by
\begin{equation*}
\begin{array}{llll}
u=\frac{U}{\mathcal E}, &w=\frac{\mathcal K\mathcal B}{\mathcal E^2\mathcal A}W, & x=\sqrt{\frac{\mathcal R} {\epsilon_w}}X, &\quad  t=\mathcal R \tau,  \\
\gamma=\frac{\mathcal K}{\mathcal E}, &\delta= \frac{
\mathcal K\mathcal R} {\mathcal E \mathcal A}, &\quad \beta=\frac{\mathcal A \mathcal H \mathcal E}{\mathcal B \mathcal K},&\quad \epsilon=\frac{\epsilon_u}{\epsilon_w}.
\end{array}
\end{equation*}
In these new variables the system \eqref{e:1h} reads
\begin{eqnarray}\label{e:2h}
u_{t}&=&\epsilon u_{xx} +\frac{1}{\delta}\left(u(\gamma-u)-\frac{uw}{1+u}\right),\notag\\
w_{t}&=& w_{xx}+w \left(1-\frac{\beta w}{u}\right).
\end{eqnarray} 
For the computational brevity we will set $\gamma=1$. 
We are interested in front solutions asymptotically connecting $(u,w)=(1,0)$ and $(u,w)=\left(\frac{-1+\sqrt{4\beta^2+1}}{2\beta}, \frac{-1+\sqrt{4\beta^2+1}}{2\beta^2} \right)$.
Switching to the moving frame $z=x-ct$ in \eqref{e:2h}, we obtain
\begin{eqnarray}\label{e:3h}
u_{t}&=&\epsilon u_{zz} +cu_z+\frac{1}{\delta}\left(u(1-u)-\frac{uw}{1+u}\right),\notag\\
w_{t}&=& w_{zz}+cw_z+w \left(1-\frac{\beta w}{u}\right).
\end{eqnarray} 
The fronts  in the system \eqref{e:2h}  are associated with heteroclinic orbits of the dynamical system
\begin{eqnarray}\label{e:4h}
u_1^{\prime}&=&u_2,\notag\\
\epsilon u_{2}^{\prime}&=&- cu_2-\frac{1}{\delta}\left(u_1(1-u_1)-\frac{u_1w_1}{1+u_1}\right),\notag\\
w_1^{\prime}&=& w_2,\notag\\
w_{2}^{\prime} &=&-cw_2-w_1 \left(1-\frac{\beta w_1}{u_1}\right),
\end{eqnarray} 
where $u_1=u$, $w_1=w$, and the derivative is taken with respect to $z$.
Assuming that  $0<\epsilon \ll \delta \ll 1$, and following the same multiple scale reduction as in Section \ref{sec:1d}, we obtain that the limiting behavior  is defined by the dynamical system 
\begin{eqnarray}\label{e:5h}
w_1^{\prime}&=& w_2,\notag\\
w_{2}^{\prime} &=&-cw_2-w_1\frac{\sqrt{1-w_1} -\beta w_1}{\sqrt{1-w_1} },
\end{eqnarray} 
which is related to the traveling wave equation in a scalar Fisher-KPP equation
\begin{eqnarray}\label{e:5hkpp}
w_t=w^{\prime\prime}+cw^{\prime}+ w\frac{\sqrt{1-w} -\beta w}{\sqrt{1-w} }
\end{eqnarray}
% with the nonlinearity $f(w)=$. 
Thus, a  theorem  holds that is similar to Theorem~\ref{T:sdp} for  the slowly diffusing  prey case for Rosenzweig-MacArthur  system:
\begin{Theorem} \label{T:sdpHT} For every fixed  $\beta>0$  and $c > 2$,  there exists  $\delta_0=\delta_0(\beta) >0$    such that for every $0<\delta <\delta_0$  there exists $\epsilon_0(\beta)>0$  such that for each $0<\epsilon <\epsilon_0$  there is a  heteroclinic orbit  of system \eqref{e:4h} 
that  asymptotically connects     the saddle  $(u_1,u_2,w_1,w_2)=(\frac{-1+\sqrt{4\beta^2+1}}{2\beta}, 0,\frac{-1+\sqrt{4\beta^2+1}}{2\beta^2},0)$  at $-\infty$  to  the saddle  $(u_1,u_2,w_1,w_2)=(1,0,0,0)$ at $+\infty$.
\end{Theorem}
The latter implies an analogue of Theorem~\ref{T:1}  for the diffusive Holling-Tanner model \eqref{e:1h}:
\begin{Theorem} \label{T:20}
For every fixed  $\beta>0$  and $c > 2$,  there exists  $\delta_0=\delta_0(\beta) >0$    such that for every $0<\delta <\delta_0$  there exists $\epsilon_0(\beta)>0$  such that for each $0<\epsilon <\epsilon_0$  there is a  translationally invariant family of  fronts  of the system \eqref{e:1h}  moving with speed $c$
that  asymptotically connect     the  equilibrium   $(u,w)=(\frac{-1+\sqrt{4\beta^2+1}}{2\beta}, \frac{-1+\sqrt{4\beta^2+1}}{2\beta^2})$  at $-\infty$  to   the equilibrium $(u,w)=(1,0)$ at $+\infty$. 
\end{Theorem}
For the diffusive Holling-Tanner model \eqref{e:1h} the case of vanishing diffusion limit when both $\epsilon_u$ and $\epsilon_w$ are very small  is considered in \cite{GMS} and the existence of fronts is  analytically proved. As in the case of the vanishing diffusion limit for Rosenzweig-MacArthur system, the connection of these fronts with Fisher-KPP fronts is not directly observed, as opposed to the regime described above.

\section{Acknowledgements and other remarks.}
During the work on this project, Ghazaryan  was supported by the NSF  grant  DMS-1311313, and  Manukian was supported  by the   Simons Foundation
through the Collaboration grant \#246535. Hong Cai worked on this project as a student at Miami University and later  as a student at Brown University.

The authors thank the anonymous reviewers  for the detailed  feedback and constructive  recommendations.

\end{document}